\documentclass[final,3p,times,num]{elsarticle}

\usepackage{amssymb}
\usepackage[english]{babel}
\usepackage{amsmath}
\usepackage{lipsum}
\usepackage{hyperref}
\usepackage{dsfont}
\usepackage{caption,subcaption}
\usepackage{physics}

\newcommand*{\logten}{\mathop{\log_{10}}}
\newcommand*{\lop}{\mathcal{L}_{x,t}}
\newcommand*{\nlop}{\mathcal{N}_{x,t}}
\newcommand*{\rop}{\mathcal{R}_{x,t}}

\newcommand{\me}{\mathrm{e}}

\journal{Applied Mathematics Letters}

\begin{document}

\begin{frontmatter}
\title{Convergence acceleration for the BLUES function method}

\author[add1]{Jonas Berx\corref{cor1}}
\ead{jonas.berx@kuleuven.be}

\address[add1]{Institute for Theoretical Physics, KU Leuven, B-3001 Leuven, Belgium}
\cortext[cor1]{Corresponding author}

\begin{abstract}
A detailed comparison is made between four different iterative procedures: Picard, Ishikawa, Mann and Picard-Krasnoselskii, within the framework of the BLUES function method and the variational iteration method. The resulting modified methods are subsequently applied to a nonlinear reaction-diffusion-advection differential equation to generate approximations to the known exact solution. The differences between the BLUES function method and the variational iteration method are illustrated by studying the approximants and the error between the obtained approximants and the exact solution.

\end{abstract}

\begin{keyword}
BLUES function method \sep variational iteration method \sep Picard iteration \sep Ishikawa iteration \sep Mann iteration \sep convergence acceleration



\end{keyword}

\end{frontmatter}

\section{Introduction}
\label{sec:intro}
Finding solutions of nonlinear differential equations in mathematics, physics and other sciences is a challenge. Because of the inherent nonlinearity, regular solution methods such as superposition or direct integration are often impossible and hence no general solution strategy exists. Therefore, it is appropriate to look for approximate solutions instead of trying to solve the differential equations directly. To this end, (semi-)analytical iterative methods such as the variational iteration method (VIM) \cite{HE1999699,HE20073,KHURI201428}, Adomian decomposition method (ADM) \cite{adomian,adomian1986nonlinear}, homotopy perturbation method (HPM) \cite{HE1999257} or BLUES function method \cite{BLUES,Berx_2020,Berx2021_PDE,Berx2021_SIRS} have been proposed. 

The standard formulations of these methods all rely on iterating through a sequence of approximants by direct substitution of the $n$th approximant into some nonlinear operator to generate the $(n+1)$th approximant. While this so-called Picard iteration is straightforward and provides a direct way to generate solutions, its convergence is often slow and its accuracy insufficient to provide useful approximations with limited computational resources.

In this paper we propose to reformulate the BLUES function method and the VIM to incorporate three other iterative procedures: Mann's procedure, Ishikawa's procedure and a hybrid Picard-Krasnoselskii procedure. Subsequently, we compare the approximants generated by each of these methods and study whether the different iterative schemes result in an increased accuracy with respect to Picard iteration for a fixed number of iterations.

The setup of this paper is as follows. In section \ref{sec:pde} we introduce a nonlinear reaction-diffusion-advection PDE along with the approximate solutions obtained by regular Picard iteration for the BLUES function method abd the VIM. Next, in section \ref{sec:iterative_processes}, we set up the aforementioned iterative procedures in a general formulation. In section \ref{subsec:other_procedures_ramos}, the different procedures are applied to the reaction-diffusion-advection PDE for both the BLUES function method and the VIM. A comparison is made between the methods as well as between the different iterative procedures, and the results are presented by studying the error between the approximants and the exact solution. Finally, in section \ref{sec:conclusions}, we present the conclusions and a future outlook.

\section{Nonlinear reaction-diffusion-advection equation}
\label{sec:pde}
We set the stage by considering a nonlinear reaction-diffusion-advection PDE \cite{ramos,Berx2021_PDE} that can be used to describe, e.g., the propagation of a temperature or chemical concentration $u$ through the combined mechanisms of diffusion, nonlinear advection and reaction, 
\begin{equation}
    \label{eq:ramos_pde}
    \begin{split}
    \nlop u &\equiv u_t - u_{xx} + u u_x +u (u+a) = 0
    \end{split}
\end{equation}
defined on $(x,t) \in \mathbb{R} \times [0,\infty)$ and $a\in\mathbb{R}$, where $\nlop$ is a nonlinear operator that is defined through its action on $u$. Subscripts such as $u_t,\, u_{xx},...$ denote differentiation with respect to the subscript variable. We consider a decaying exponential initial condition,
\begin{equation}
    \label{eq:ramos_ic}
    u(x,0) \equiv f(x) = \me^{-x}\, .
\end{equation}
This unbounded initial condition is unphysical but will serve as an ideal testbed for the comparison of the different approximation methods, as in this case a simple exact solution of \eqref{eq:ramos_pde} can be found, i.e.,
\begin{equation}
    \label{eq:ramos_exact}
    u(x,t) = \me^{-(x+(a-1)t)}\, .
\end{equation}
It can be shown \cite{Berx2021_PDE} that for regular Picard iteration the VIM and BLUES function method give the following $n$th-order approximants to the exact solution \eqref{eq:ramos_exact}, respectively
\begin{subequations}
    \label{eq:sequences}
    \begin{align}
        u^{(n)}_{VIM}(x,t) &= \me^{-x}\sum_{i=0}^n \frac{(-(a-1)t)^i}{i!}\label{eq:sequence_vim}\\
        u^{(n)}_{BLUES}(x,t) &= \me^{-at-x} \sum_{i=0}^n\frac{t^i}{i!}\, .\label{eq:sequence_blues}
    \end{align}
\end{subequations}
It is now easy to see that, in the limit for $n\rightarrow\infty$, both methods converge to the exact solution \eqref{eq:ramos_exact}.

\section{Iterative methods and procedures}
\label{sec:iterative_processes}

To avoid confusion, we will from now on differentiate between \textit{methods} and \textit{procedures}. A method, e.g. the BLUES function method or the VIM, signifies a particular strategy employed to find an approximate solution $u_n$. A procedure, however, is the specific manner in which one iterates through the approximants generated by an iterative method. Examples are Picard or Mann iteration.

Consider a sequence of functions $\{u_n\}_{n=0}^\infty$ and a (nonlinear) mapping $\mathcal{T}: C\rightarrow C$, with $C$ a nonempty convex subset of a normed linear space $E$. The Picard iterative procedure is the following
\begin{equation}
    \label{eq:methods_picard}
    \begin{split}
        u_{n+1} = \mathcal{T}[u_n]\,, \qquad u_0 = u\, \in C\,,
    \end{split}
\end{equation}
with $n \in\mathbb{N}$. 

For the BLUES function method, the action of the operator $\mathcal{T}$ on the $n$th approximant $u_n$ is
\begin{equation}
    \label{eq:methods_operators_T_blues}
    \begin{split}
        \mathcal{T}[u_n] &= u_0 + G\ast\mathcal{R}[u_n]\,,
    \end{split}
\end{equation}
where $G$ is the BLUES function, which we take to be the Green function of a judiciously chosen linear problem, i.e., $\lop G =\delta(t)$, with the same initial condition $u(x,0) = f(x)$, where $\lop$ is a linear operator. For the remainder of this work, we will assume that $\lop$ is an operator that at most contains first-order time derivatives. The $*$-multiplication indicates the convolution operator, while $\rop\equiv\lop-\nlop$ is the residual operator. A detailed setup of the BLUES function method can be found in e.g., Ref. \cite{Berx2021_PDE}.

Now, for the VIM, the action of the operator $\mathcal{T}$ can be written as
\begin{equation}
    \label{eq:methods_operators_T_vim}
    \begin{split}
        \mathcal{T}[u_n] &= u_n + \int_{t_0}^t\mathrm{d}s\,\lambda(s) (Lu_n + N\tilde{u}_n -f)\,,
    \end{split}
\end{equation}
where $L,N$ and $f$ are respectively the linear operator, which is often the highest time derivative, the nonlinear operator, and an external source. The function $\lambda(s)$ is a general Lagrange multiplier that can be identified via variational theory, and $\tilde{u}_n$ is a restricted variation, i.e., $\delta\tilde{u}=0$. For the VIM, $u_0 = u(x,0)$.

We now briefly discuss other iterative procedures that can be embedded into these two methods.

\subsection{Mann's iterative procedure}
\label{subsec:mann}

Consider the following single-step procedure
\begin{equation}
    \label{eq:methods_mann}
    \begin{split}
        u_{n+1} = (1-\alpha_n) u_n + \alpha_n \mathcal{T}[u_n]\,, \qquad u_0 = u\,\in C\,,
    \end{split}
\end{equation}
with $n \in\mathbb{N}$ and where $\left(\alpha_n\right)_{n\in\mathbb{N}}$ is a sequence of positive real numbers. This process is called Mann's iterative procedure \cite{mann1953}. This scheme is sometimes used to accelerate the convergence of the VIM by considering the $\alpha_n$'s as convergence-control parameters \cite{Ahmad2020,Hosseini2010}. These can be optimally determined by minimising the residual square error of the approximants with respect to the $\alpha_n$ in each iteration. This is numerically expensive but can result in needing a smaller number of approximants to achieve the desired accuracy. When $\alpha_n = \alpha\, \text{(constant)}\,, \forall n\in\mathbb{N}$, the procedure is called Krasnoselskii's iterative procedure, while for $\alpha_n = 1$, Mann's iteration \eqref{eq:methods_mann} reduces to Picard iteration \eqref{eq:methods_picard}.

\subsection{Ishikawa iterative procedure}
\label{subsec:ishikawa}

Consider the following two-step procedure
\begin{equation}
    \label{eq:methods_Ishikawa}
    \begin{split}
        u_{n+1} &= (1-\alpha_n) u_n + \alpha_n \mathcal{T}[v_n]\,, \qquad  u_0 = u\, \in C\\
        v_n &= (1-\beta_n) u_n + \beta_n \mathcal{T}[u_n]\,,
    \end{split}
\end{equation}
with $n \in\mathbb{N}$ and where $(\alpha_n)_{n\in\mathbb{N}}$ and $(\beta_n)_{n\in\mathbb{N}}$ are sequences of positive real numbers.  This process is called the Ishikawa iterative procedure \cite{Ishikawa1974,KHURI201850}. As was the case for Mann's procedure, the parameters $\alpha_n$ and $\beta_n$ can be used as convergence-control parameters. Obviously, for $\beta_n=0$, Ishikawa's procedure \eqref{eq:methods_Ishikawa} reduces to Mann's procedure \eqref{eq:methods_mann}.

\subsection{Hybrid Picard-Krasnoselskii procedure}
\label{subsec:hybrid}
While the choice for any one of the above procedures may improve convergence for the BLUES function method or the VIM, there exists another choice by combining Picard's and Krasnoselskii's method. This so-called hybrid Picard-Krasnoselskii's procedure can be proven to converge faster than Picard's, Mann's, Krasnoselskii's or Ishikawa's procedures \cite{Okeke2017}. It can be described as follows,
\begin{equation}
    \label{eq:methods_hybrid}
    \begin{split}
        u_{n+1} &= \mathcal{T}[v_n]\,,\qquad u_0 = u\, \in C\\
        v_n &= (1-\lambda) u_n + \lambda\, \mathcal{T}[u_n]\,,
    \end{split}
\end{equation}
with $n \in\mathbb{N}$. Again, the parameter $\lambda$ can be used to control the convergence. For $\lambda=0$, the hybrid procedure reduces to Picard iteration, while for $\lambda=1$, it becomes a two-step Picard iteration, where the $(n+1)$th approximant is found by applying the $\mathcal{T}$ operator twice. 

\section{Comparison}
\label{subsec:other_procedures_ramos}
Let us now consider the procedures that were discussed in section \ref{sec:iterative_processes}. The operator $\mathcal{T}$ from equations \eqref{eq:methods_operators_T_blues} and \eqref{eq:methods_operators_T_vim} for the nonlinear reaction-diffusion-advection equation \eqref{eq:ramos_pde} is for respectively the BLUES function method and the VIM
\begin{equation}
    \label{eq:methods_operators_T_full}
    \begin{split}
        \mathcal{T}[u^{(n)}] &= \me^{-at -x} + \int\limits_{0}^t \mathrm{d}s\, G(t-s)\left[\pdv[2]{u^{(n)}}{x} - u^{(n)}\pdv{u^{(n)}}{x} -(u^{(n)})^2\right](x,s)\,,\\
        \mathcal{T}[u_n] &= u_n - \int\limits_0^t \mathrm{d}s \left[\pdv{u_n}{s} - \pdv[2]{u_n}{x} + u_n\pdv{u_n}{x} +u_n(u_n+a)\right](x,s)\,.
    \end{split}
\end{equation}
Note that for the BLUES method we have used superscripts to denote the iterations, while for the VIM we use subscripts. We will use these interchangeably when no confusion can arise. The function $G(t) = \me^{-a(t)}\Theta(t)$, where $\Theta(.)$ is the Heaviside function, is the Green function of the associated linear operator $\lop u\equiv u_t + a u$. This choice of linear operator fixes the residual operator to be $\rop u\equiv \lop u-\nlop u = u_{xx} - u u_x - u^2$. We now discuss the different procedures first for the VIM and subsequently for the BLUES function method.
\subsection{Variational iteration method}
For the VIM, the Picard iteration approximants are given by equation \eqref{eq:sequence_vim}. Mann's iterative procedure can be set up as follows
\begin{equation}
    \label{eq:ramos_VIM_mann}
    u_{n+1}(x,t) = u_n(x,t) - \alpha_n \int\limits_0^t \mathrm{d}s \left(\pdv{u_n}{s} - \pdv[2]{u_n}{x} + u_n\pdv{u_n}{x} +u_n(u_n+a)\right)\,.
\end{equation}
Consider now the residual error $\mathcal{E}(n,x,T)$,
\begin{equation}
    \label{eq:squared_error_definition}
    \mathcal{E}(n,x,T) = \frac{1}{T}\int\limits_0^T\left[\pdv{u_n}{s} - \pdv[2]{u_n}{x} + u_n\pdv{u_n}{x} +u_n(u_n+a)\right]^2\mathrm{d}s\,.
\end{equation}
One has to minimize $\mathcal{E}(n,x,T)$ with respect to $\alpha_n$ for each $n$. With the (arbitrary) choice $T=1$, the $\alpha_n$ with $n\in\{1,2,3\}$ become, for $a=2$,
\begin{equation}
    \label{eq:ramos_VIM_mann_alphas}
    \begin{split}
        \alpha_1 = 0.64286\,, \qquad \alpha_2 = 0.84990\,, \qquad \alpha_3 = 0.74308\,.
    \end{split}
\end{equation}
Note that the coordinate $x$ does not have to be fixed for the minimisation of the residual error. By careful inspection of the different iterative procedures \eqref{eq:methods_picard}, \eqref{eq:methods_mann}, \eqref{eq:methods_Ishikawa} and \eqref{eq:methods_hybrid}, one can see that the exponential containing the $x-$coordinate can be divided out when minimising \eqref{eq:squared_error_definition}.

For the Ishikawa scheme, we set up the procedure as follows
\begin{equation}
    \label{eq:ramos_VIM_ishikawa}
    \begin{split}
        v_n(x,t) &= u_n(x,t) - \beta_n \int\limits_0^t \mathrm{d}s \left(\pdv{u_n}{s} - \pdv[2]{u_n}{x} + u_n\pdv{u_n}{x} +u_n(u_n+a)\right)\\
        u_{n+1}(x,t) &= (1-\alpha_n) u_n(x,t) + \alpha_n v_n(x,t)+ \alpha_n \int\limits_0^t \mathrm{d}s \left(\pdv{v_n}{s} - \pdv[2]{v_n}{x} + v_n\pdv{v_n}{x} +v_n(v_n+a)\right)
    \end{split}
\end{equation}

Minimizing the residual error $\mathcal{E}(n,x,T)$ with respect to $\alpha_n$ and $\beta_n$, with once again $T=1$ and $n\in\{1,2,3\}$ gives the following convergence-control parameters, for $a=2$,
    \begin{align}
        \alpha_1 &= 0.94272\,, & \alpha_2 &= 0.97810\,, & \alpha_3 &= 0.98356\,,\label{eq:ramos_VIM_ishikawa_alphass}\\
        \beta_1 &= 0.65972\,, & \beta_2 &= 0.84372\,, & \beta_3 &= 0.76055\,.\label{eq:ramos_VIM_ishikawa_betas}
    \end{align}
Note that the $\alpha_n$ converge to unity, indicating that the Ishikawa procedure converges to a hybrid Picard-Ishikawa procedure. Lastly, the Picard-Krasnoselskii hybrid procedure is the following
\begin{equation}
    \label{eq:ramos_VIM_PicKras}
    \begin{split}
        v_n(x,t) &= u_n(x,t) - \lambda\int\limits_0^t \mathrm{d}s \left(\pdv{u_n}{s} - \pdv[2]{u_n}{x} + u_n\pdv{u_n}{x} +u_n(u_n+a)\right)\,,\\
        u_{n+1}(x,t) &= v_n(x,t) - \int\limits_0^t \mathrm{d}s \left(\pdv{v_n}{s} - \pdv[2]{v_n}{x} + v_n\pdv{v_n}{x} +v_n(v_n+a)\right)\,.\\
    \end{split}
\end{equation}
For $n=3$, the parameter $\lambda$ can be calculated by minimising once again the residual error with respect to $\lambda$, resulting in $\lambda = 0.85590$ for $T=1$ and $a=2$.

\begin{figure}[htp]
    \centering
    \begin{subfigure}{0.465\textwidth}
        \includegraphics[width=\linewidth]{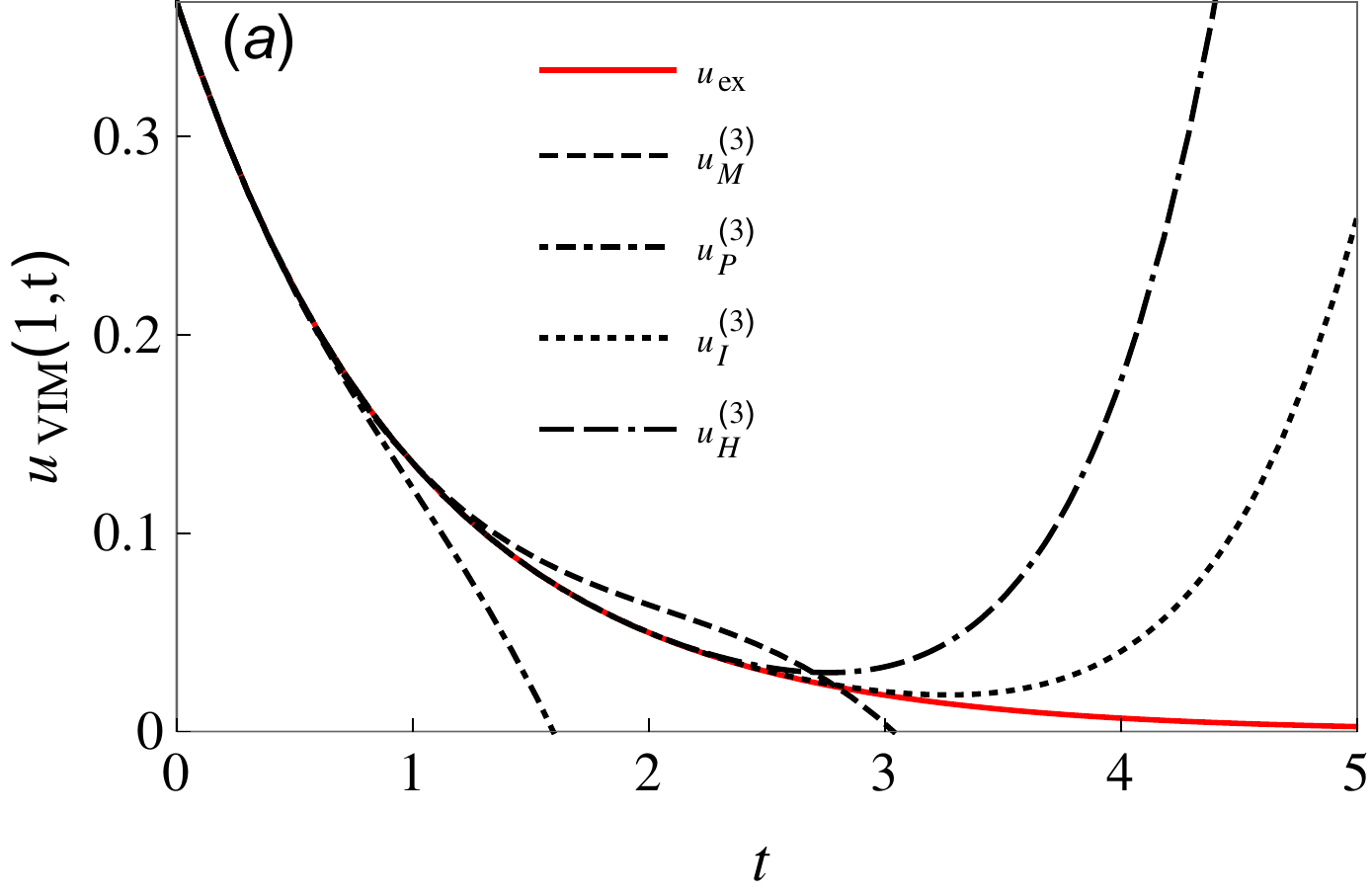}
    \end{subfigure}
    \begin{subfigure}{0.495\textwidth}
        \includegraphics[width=\linewidth]{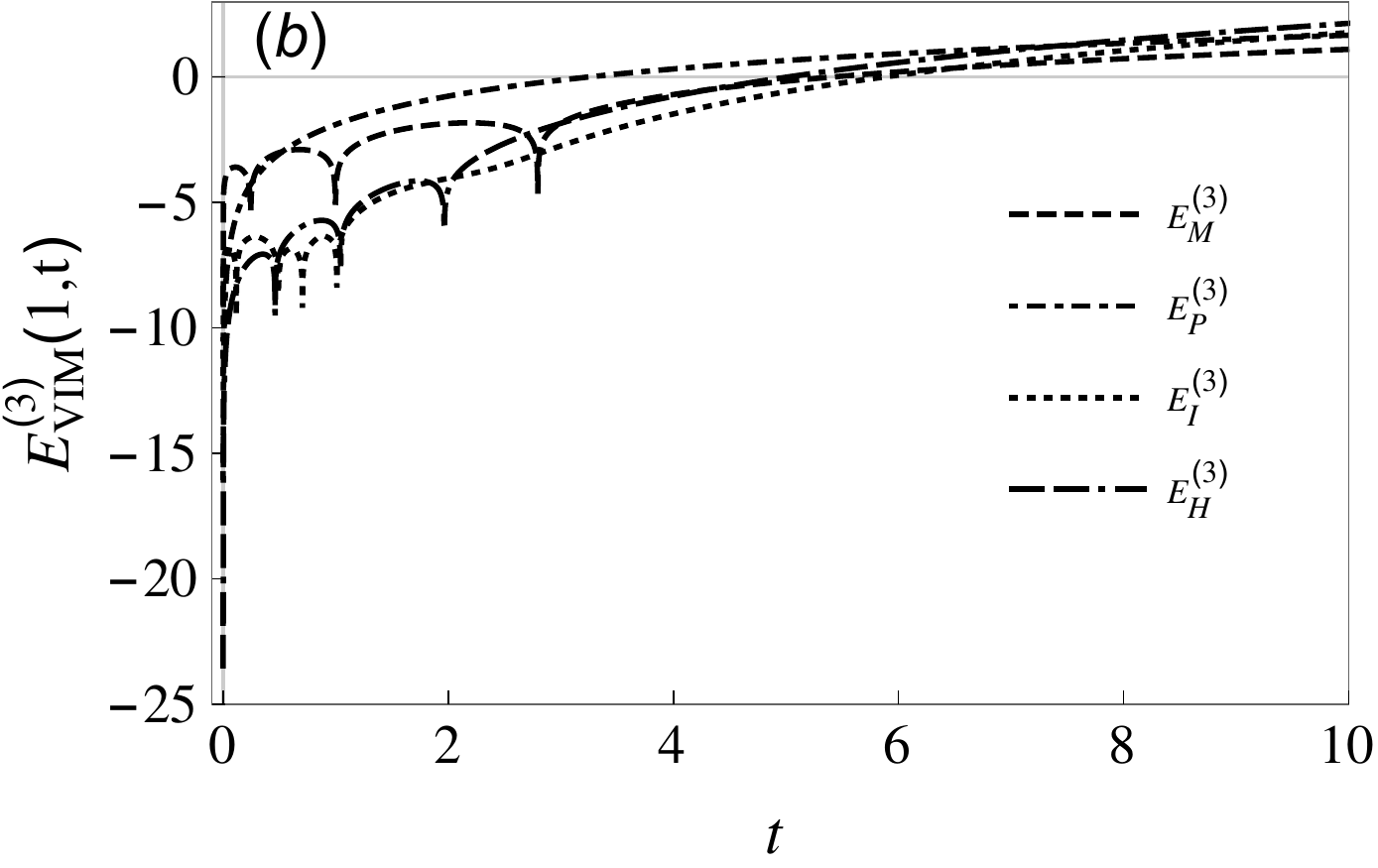}
    \end{subfigure}
    \caption{\textbf{(a)} Comparison of the different iteration procedures applied to the third-order VIM, with $a=2$. Mann's procedure \eqref{eq:ramos_VIM_mann} (dashed line), Picard's procedure \eqref{eq:sequence_vim} (dot-dashed line), Ishikawa's procedure \eqref{eq:ramos_VIM_ishikawa} (dotted line), the Picard-Krasnoselskii hybrid scheme \eqref{eq:ramos_VIM_PicKras} (dot-dash-dashed line) and the exact solution \eqref{eq:ramos_exact} (red line). \textbf{(b)} Comparison of the third-order VIM errors $E^{(3)}(x=1,t)$. The spatial coordinate is fixed at $x=1$.}
    \label{fig:ramos_methods_VIM_solution}
\end{figure}

The (logarithmic) error of the different procedures with respect to the exact solution can be defined as
\begin{equation}
    \label{eq:ramos_VIM_errors}
    E^{(n)}(x,t) = \logten \left|u_{ex}(x,t)-u^{(n)}(x,t)\right|\,,
\end{equation}
and will be used to study the differences between the procedures. The different iteration procedures applied to the VIM for the solution of the nonlinear reaction-diffusion-convection equation \eqref{eq:ramos_pde} are compared in Fig. \ref{fig:ramos_methods_VIM_solution}(a), where the errors between the exact solution and the approximants for the different procedures are shown in panel (b). 

\subsection{BLUES function method}
Let us now study the application of the different iterative procedures on the BLUES function method. The Picard iteration approximants are given by equation \eqref{eq:sequence_blues}. Mann's iterative procedure can be set up as follows
\begin{equation}
    \label{eq:ramos_BLUES_mann}
    \begin{split}
    u^{(n+1)}(x,t) &= (1-\alpha_n)u^{(n)}(x,t) + \alpha_n u^{(0)}(x,t)+ \alpha_n \int\limits_{0^-}^t \mathrm{d}s\, G(t-s)\left[u^{(n)}_{xx} - u^{(n)}u^{(n)}_x -(u^{(n)})^2\right]
    \end{split}
\end{equation}

We now repeat the same strategy as before and minimize $\mathcal{E}(n,x,T)$ with respect to $\alpha_n$ for each $n$. With the choice $T=1$, the $\alpha_n$ with $n\in\{1,2,3\}$ become, for $a=2$,
\begin{equation}
    \label{eq:ramos_BLUES_mann_alphas}
    \begin{split}
        \alpha_1 = 1.2118\,,\qquad
        \alpha_2 = 1.1993\,,\qquad
        \alpha_3 = 1.0524\,.
    \end{split}
\end{equation}
The Ishikawa scheme for the BLUES function method can be set up as follows
\begin{equation}
    \label{eq:ramos_BLUES_ishikawa}
    \begin{split}
        v^{(n)}(x,t) &= (1-\beta_n) u^{(n)}(x,t) + \beta_n u^{(0)}(x,t)
        + \beta_n\int\limits_{0^-}^t \mathrm{d}s\, G(t-s)\left[u^{(n)}_{xx} - u^{(n)}u^{(n)}_x -(u^{(n)})^2\right]\,,\\
        u^{(n+1)}(x,t) &= (1-\alpha_n) u^{(n)}(x,t) + \alpha_n u^{(0)}(x,t)
        + \alpha_n \int\limits_{0^-}^t \mathrm{d}s\, G(t-s)\left[v^{(n)}_{xx} - v^{(n)}v^{(n)}_x -(v^{(n)})^2\right]\,.
    \end{split}
\end{equation}
The convergence-control parameters are, for $a=2$,
\begin{align}
    \label{eq:ramos_BLUES_ishikawa_alphas_betas}
        \alpha_1 &= 0.95312\,, & \alpha_2 &= 0.98052\,, & \alpha_3 &= 0.99276\,,\\
        \beta_1 &= 1.52507\,, & \beta_2 &= 1.14221\,, & \beta_3 &= 1.1600\,.
\end{align}
Note that once again the $\alpha_n$ converge to unity, so this Ishikawa procedure for the BLUES function method also converges to a hybrid Picard-Ishikawa procedure. 

Next, the Picard-Krasnoselskii hybrid procedure is the following
\begin{equation}
    \label{eq:ramos_BLUES_PicKras}
    \begin{split}
        v^{(n)}(x,t) &= (1-\lambda) u^{(n)}(x,t) + \lambda u^{(0)}(x,t) 
        + \lambda\int\limits_{0^-}^t \mathrm{d}s\, G(t-s)\left[u^{(n)}_{xx} - u^{(n)}u^{(n)}_x -(u^{(n)})^2\right]\,,\\
        u^{(n+1)}(x,t) &= u^{(0)}(x,t) + \int\limits_{0^-}^t \mathrm{d}s\, G(t-s)\left[v^{(n)}_{xx} - v^{(n)}v^{(n)}_x -(v^{(n)})^2\right]\,.
    \end{split}
\end{equation}
For $n=3$, the parameter $\lambda$ is found as $\lambda = 1.18323$ for $T=1$ and $a=2$.

\begin{figure}[htp]
    \centering
    \begin{subfigure}{0.465\linewidth}
        \includegraphics[width=\linewidth]{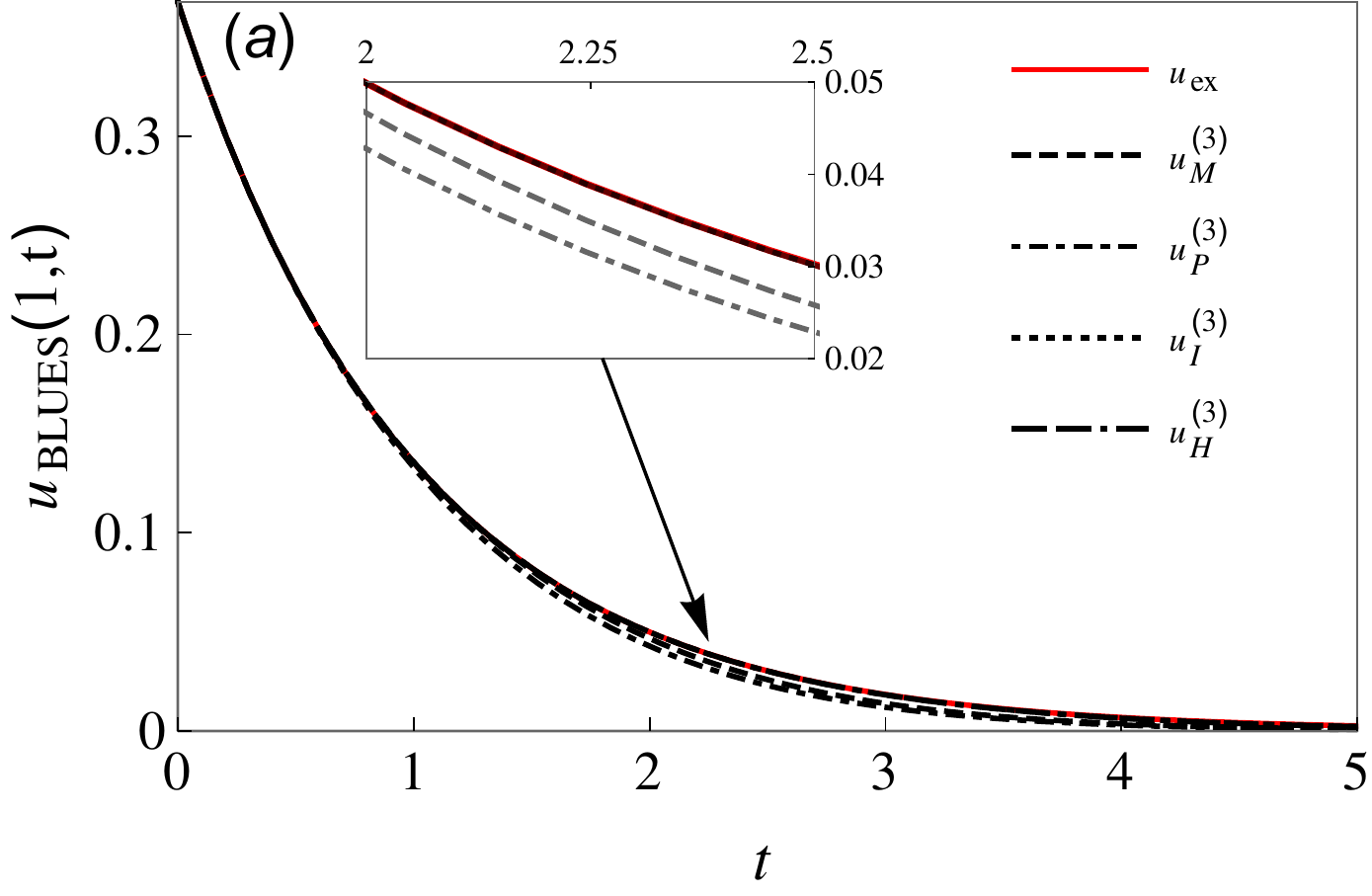}
    \end{subfigure}
    \begin{subfigure}{0.495\linewidth}
        \includegraphics[width=\linewidth]{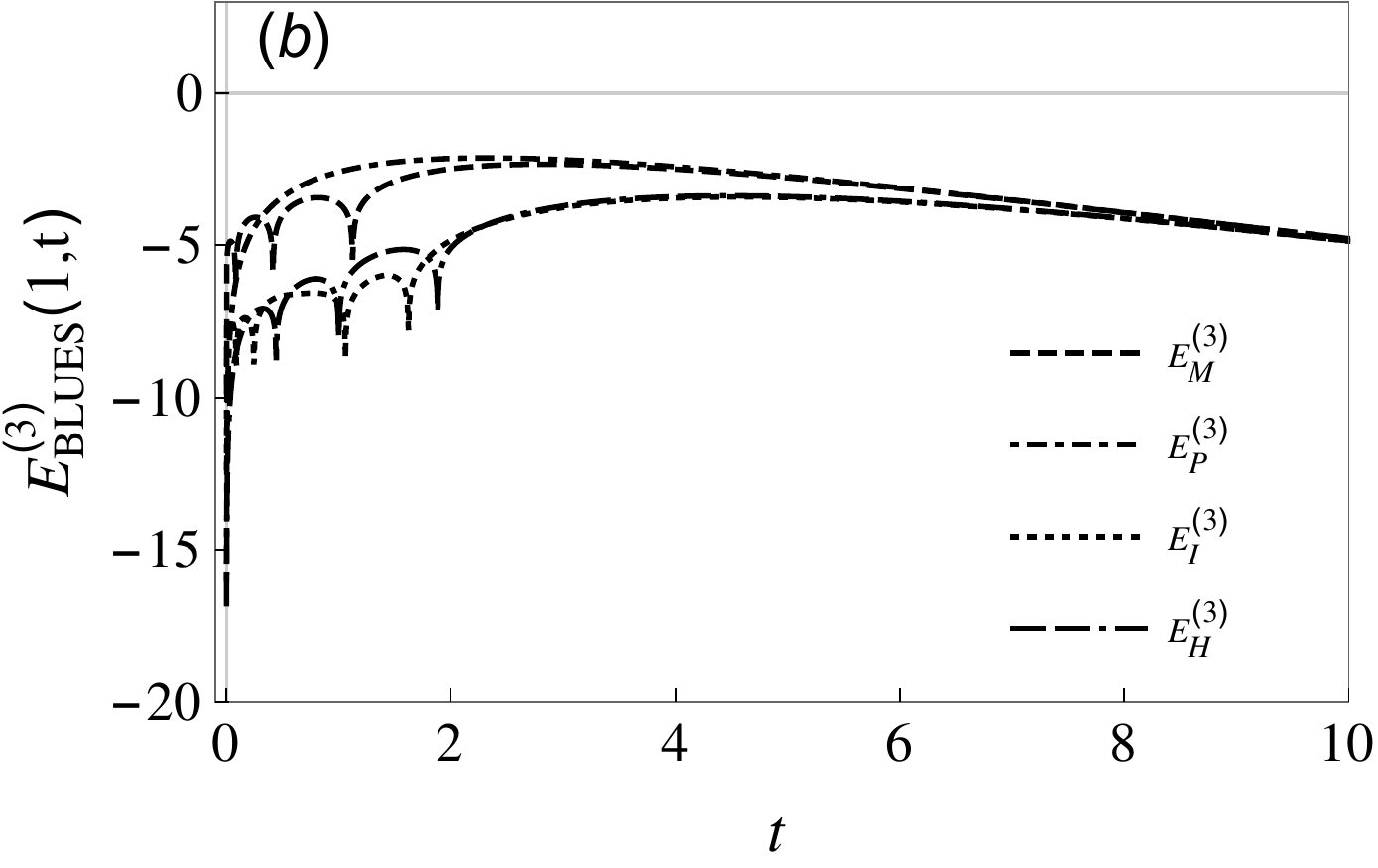}
    \end{subfigure}
    \caption{\textbf{(a)} Comparison of the different iteration procedures applied to the third-order BLUES function method, with $a=2$. Mann's procedure \eqref{eq:ramos_BLUES_mann} (dashed line), Picard's procedure \eqref{eq:sequence_blues} (dot-dashed line), Ishikawa's procedure \eqref{eq:ramos_BLUES_ishikawa} (dotted line), the Picard-Krasnoselskii hybrid scheme \eqref{eq:ramos_BLUES_PicKras} (dot-dash-dashed line) and the exact solution \eqref{eq:ramos_exact} (red line). A zoomed-in view shows that both Ishikawa's procedure and the hybrid scheme coincide almost perfectly with the exact solution. \textbf{(b)} Comparison of the third-order BLUES function method errors $E^{(3)}(x=1,t)$. The spatial coordinate is fixed at $x=1$.}
    \label{fig:ramos_methods_BLUES_solution}
\end{figure}

The different iteration procedures applied to the BLUES function method for the solution of the nonlinear reaction-diffusion-convection equation \eqref{eq:ramos_pde} are compared in Fig. \ref{fig:ramos_methods_BLUES_solution}. 

It is clear from Figure \ref{fig:ramos_methods_BLUES_solution} that Mann's and Picard's schemes converge to the same iterative procedure. Moreover, this is also the case for Ishikawa's scheme and the hybrid Picard-Krasnoselskii scheme. For $t\rightarrow\infty$, all methods converge and the error decreases exponentially. Comparing the results obtained by the VIM and by the BLUES function method, it becomes clear that BLUES outperforms the VIM for longer times $t$.

\section{Conclusions}
\label{sec:conclusions}
In this work, we have applied four different iteration procedures (Picard, Mann, Ishikawa and hybrid Picard-Krasnoselskii) to both the variational iteration method and the BLUES function method in the framework of a nonlinear reaction-diffusion-advection partial differential equation. By numerically determining the optimal convergence-control parameters, we have shown that for the aforementioned equation, the Ishikawa and hybrid Picard-Krasnoselskii schemes produce significantly more accurate approximations to the exact solution, for an equal number of iterations. Moreover, for the BLUES function method, the different schemes all produce globally convergent approximants, for which the difference between the approximants and the exact solution decreases exponentially in time. While implementing other iterative procedures into the VIM significantly improves its convergence and reduces the number of iterations required, the BLUES method seems to outperform the VIM for longer times. 

We envision implementing these iterative methods into more general formulations of the BLUES function method such as e.g., higher-order PDEs, coupled PDEs, or stochastic differential equations (SDEs). 

\section*{Acknowledgement}
The author wishes to express his gratitude to his colleague Joseph O. Indekeu, without whom this research would not have been possible.

This research did not receive any specific grant from funding agencies in the public, commercial, or not-for-profit sectors.

\bibliographystyle{elsarticle-num}
\bibliography{BLUES.bib}

\end{document}